\documentclass[12pt]{article}
\usepackage{amsmath}
\usepackage{amssymb}
\usepackage{vatola}

\def\q{\quad}
\def\qq{\qquad}
\def\qtq#1{\q\t{#1}\q}
 \TagsOnRight
\def\t{\hbox}
\def\mod#1{\ (\text{\rm mod}\ #1)}
\def\f{\frac}
\def\e{\equiv}
\def\b{\binom}

\begin{document}
\par Bull. Aust. Math.
Soc. 86(2012), no.1, 164-176
\par\q
\let \pro=\proclaim
\let \endpro=\endproclaim

\centerline {\Large\bf
 Ramsey numbers for trees\q}
$$\q$$
\centerline{\bf \qq Zhi-Hong Sun} \centerline{School of Mathematical
Sciences,} \centerline{Huaiyin Normal University,}
\centerline{Huaian, Jiangsu 223001, P.R. China}
 \centerline{E-mail:
zhihongsun@yahoo.com} \centerline{Homepage:
http://www.hytc.edu.cn/xsjl/szh}
\abstract{For $n\ge 5$ let
  $T_n'$ denote the unique tree on $n$ vertices
with $\Delta(T_n')=n-2$, and let $T_n^*=(V,E)$ be the tree on $n$
vertices with
 $V=\{v_0,v_1,\ldots,$ $v_{n-1}\}$ and
 $E=\{v_0v_1,\ldots,v_0v_{n-3},$ $v_{n-3}v_{n-2},v_{n-2}v_{n-1}\}$. In
 this paper we evaluate the Ramsey numbers $r(G_m,T_n')$ and
 $r(G_m,T_n^*)$, where $G_m$ is a connected graph of order $m$.
  As examples, for $n\ge 8$ we have
  $r(T_n',T_n^*)=r(T_n^*,T_n^*)=2n-5$, for $n>m\ge 7$
    we have $r(K_{1,m-1},T_n^*)=m+n-3$ or $m+n-4$ according as
    $m-1\mid (n-3)$ or $m-1\nmid (n-3)$,
 for $m\ge 7$ and $n\ge (m-3)^2+2$ we have $r(T_m^*,T_n^*)=m+n-3$ or
 $m+n-4$ according as $m-1\mid (n-3)$ or $m-1\nmid (n-3)$.
\par\q
 }
 \endabstract
  \footnotetext[1] {The author is
supported by the National Natural Sciences Foundation of China
(grant no. 10971078).}

\section*{1. Introduction}\par\q
\par In this paper, all graphs are simple graphs. For a graph $G=(V(G),E(G))$
let $e(G)=|E(G)|$ be the number of edges in $G$ and let
$\Delta(G)$ be the maximal degree of $G$.
  For a forbidden graph $L$, let $ex(p;L)$
denote the maximal number of edges in a graph of order $p$ not
containing $L$ as a subgraph. The corresponding Tur\'an's problem is
to evaluate $ex(p;L)$.
\par Let $\Bbb N$ be the set of
positive integers, and let $p,n\in\Bbb N$ with $p\ge n\ge 3$. For a
given tree $T_n$ on $n$ vertices,
 it is difficult to determine the value of $ex(p;T_n)$. The famous Erd\"os-S\'os
 conjecture asserts that
$ex(p;T_n)\le \f{(n-2)p}2$ for every tree $T_n$ on $n$ vertices. For
the progress on the Erd\"os-S\'os
 conjecture, see [4,8,9,11].
 Write
$p=k(n-1)+r$, where $k\in\Bbb N$ and $r\in\{0,1,\ldots,n-2\}$. Let
$P_n$ be the path on $n$ vertices. In [5] Faudree and Schelp showed
that
$$ex(p;P_n)=k\binom {n-1}2+\binom r2=\f{(n-2)p-r(n-1-r)}2.\tag 1.1$$
In the special case  $r=0$, (1.1) is due to Erd$\ddot{\t{\rm o}}$s
and Gallai [3]. Let $K_{1,n-1}$ denote the unique tree on $n$
vertices with $\Delta(K_{1,n-1})=n-1$, and for $n\ge 4$ let
  $T_n'$ denote the unique tree on $n$ vertices
with $\Delta(T_n')=n-2$. In [10] the author and Lin-Lin Wang
obtained exact values of $ex(p;K_{1,n-1})$ and $ex(p;T_n')$, see
Lemmas 2.4 and 2.5.
\par For $n\ge 5$ let $T_n^*=(V,E)$ be the tree on $n$ vertices with
 $V=\{v_0,v_1,\ldots,v_{n-1}\}$ and
 $E=\{v_0v_1,\ldots,v_0v_{n-3},v_{n-3}v_{n-2},v_{n-2}v_{n-1}\}$.
In [10], we also determine the value of $ex(p;T_n^*)$, see Lemmas
2.6-2.8.
\par  As usual $\overline{G}$ denotes the
complement of a graph $G$. Let $G\sb 1$ and $G\sb 2$ be two graphs.
The Ramsey number $r(G\sb 1, G\sb 2)$ is the smallest positive
integer $n$ such that, for every graph $G$ with $n$ vertices, either
$G$ contains  a copy of $G\sb 1$
 or else $\overline{G}$ contains a copy of $G_2$.
 \par Let $n\in\Bbb N$ with $n\ge 6$. If the Erd\"os-S\'os
 conjecture is true, it is known that $r(T_n,T_n)\le 2n-2$ (see [8]).
 Let $m,n\in\Bbb N$. In 1973 Burr and
 Roberts[2] showed that for $m,n\ge 3$,
$$r(K_{1,m-1},K_{1,n-1})=\cases m+n-3&\t{if $2\nmid mn$,}
\\m+n-2&\t{if $2\mid mn$}.\endcases$$
In 1995, Guo and Volkmann[6] proved that for $n\ge m\ge 5$,
$$r(T_m',T_n')=\cases m+n-3&\t{if $m-1\mid (n-3)$,}
\\m+n-5&\t{if $m=n\e 0\mod 2$,}
\\m+n-4&\t{otherwise}\endcases$$
and, for $n>m\ge 4$,
$$r(K_{1,m-1},T_n')=\cases m+n-3&\t{if $2\mid m(n-1)$,}
\\m+n-4&\t{if $2\nmid m(n-1)$.}
\endcases$$
\par Let $m,n\in\Bbb N$ with $n\ge m\ge 6$.
In this paper we evaluate the Ramsey number $r(T_m,T_n^*)$ for
 $T_m\in\{P_m,K_{1,m-1},$ $T_m',T_m^*\}$. As examples,
for $n\ge 8$,
  $$r(P_n,T_n^*)=r(T_n^*,T_n^*)=2n-5;$$ for $n>m\ge 7$,
    $$r(K_{1,m-1},T_n^*)=\cases m+n-3&\t{if $m-1\mid (n-3)$,}
\\m+n-4&\t{if $m-1\nmid (n-3)$;}\endcases$$
 and, for $m\ge 7$ and $n\ge (m-3)^2+2$,
 $$r(P_m,T_n^*)=r(T_m',T_n^*)=r(T_m^*,T_n^*)=
\cases m+n-3&\t{if $m-1\mid (n-3)$,}
\\m+n-4&\t{if $m-1\nmid (n-3)$.}\endcases$$

 \par In addition to the above notation, throughout the paper we
also use the following notation:
  $\lfloor x\rfloor$ is the greatest integer not
exceeding $x$, $K_n$ is the complete graph on $n$ vertices,
$K_{m,n}$ is the complete bipartite graph with $m$ and $n$ vertices
in the bipartition, $d_G(v)$ is the degree of the vertex $v$ in
given graph $G$,  and $d(u,v)$ is the distance between the two
vertices $u$ and $v$ in a graph.

 \section*{2. Basic lemmas}\par\q
\pro{Lemma 2.1} Let $G_1$ and $G_2$ be two graphs. Suppose $p\in\Bbb
N, p\ge max\{|V(G_1)|$, $|V(G_2)|\}$ and $ex(p;G_1)+ex(p;G_2)<\b
p2.$ Then $r(G_1,G_2)\le p.$\endpro Proof. Let $G$ be a graph of
order $p$. If $e(G)\le ex(p;G_1)$ and $e(\overline{G})\le
ex(p;G_2)$, then
$$ex(p;G_1)+ex(p;G_2)\ge e(G)+e(\overline{G})=\b p2.$$ This contradicts
the assumption. Hence, either $e(G)>ex(p;G_1)$ or
$e(\overline{G})>ex(p;G_2).$ Therefore, $G$ contains a copy of $G_1$
or $\overline{G}$ contains a copy of $G_2$. This shows that
$r(G_1,G_2)\le |V(G)|=p.$ So the lemma is proved.

\pro{Lemma 2.2} Let $k,p\in\Bbb N$ with $p\ge k+1$. Then there
exists a $k-$regular graph of order $p$ if and only if $2\mid kp$.
\endpro
\par This is a known result; see, for example, [10,
Corollary 2.1].
 \pro{Lemma 2.3} Let $G_1$ and $G_2$ be two
graphs with $\Delta(G_1)=d_1\ge 2$ and $\Delta(G_2)=d_2\ge 2$. Then
\par $(\t{\rm i})$ $r(G_1,G_2)\ge d_1+d_2-(1-(-1)^{(d_1-1)(d_2-1)})/2$.
\par $(\t{\rm ii})$ Suppose that $G_1$ is a connected graph of order
$m$ and $d_1<d_2\le m$. Then $r(G_1,G_2)\ge 2d_2-1\ge d_1+d_2$.
\par $(\t{\rm iii})$ Suppose that $G_1$ is a connected graph of order
$m$ and $d_2>m$.  If one of the conditions
\par \qq$(\t{\rm 1})$ $2\mid (d_1+d_2-m)$,
\par \qq$(\t{\rm 2})$ $d_1\not=m-1$,
\par \qq$(\t{\rm 3})$ $G_2$ has two
vertices $u$ and $v$ such that $d(v)=\Delta(G_2)$ and $d(u,v)=3$
\newline holds, then $r(G_1,G_2)\ge d_1+d_2$.
\endpro
Proof. We first consider (i). If $2\mid (d_1-1)(d_2-1)$, then $2\mid
(d_1-1)(d_1+d_2-1)$. Since $d_1-1\ge 1$, by Lemma 2.2 we may
construct a $d_1-1$-regular graph $G$ of order $d_1+d_2-1$. Since
$\Delta(G)=d_1-1$ and $\Delta(\overline  G)=d_2-1$, $G$ does not
contain $G_1$ as a subgraph and $\overline G$ does not contain $G_2$
as a subgraph. Hence $r(G_1,G_2)\ge 1+|V(G)|=d_1+d_2$. Now we assume
$2\nmid (d_1-1)(d_2-1)$. Then $2\mid d_1$, $2\mid d_2$ and so $2\mid
(d_1+d_2-2)$. By Lemma 2.2, we may construct a $d_1-1$-regular graph
$G$ of order $d_1+d_2-2$. Since $\Delta(G)=d_1-1$ and
$\Delta(\overline  G)=d_2-2$, $G$ does not contain $G_1$ as a
subgraph and $\overline G$ does not contain $G_2$ as a subgraph.
Hence $r(G_1,G_2)\ge 1+|V(G)|=d_1+d_2-1$. This proves (i).
\par Next we consider (ii).
Suppose that $G_1$ is a connected graph of order $m$ and $d_1<d_2\le
m$. Since $K_{d_2-1}\cup K_{d_2-1}$ does not contain any copies of
$G_1$, and its complement $K_{d_2-1,d_2-1}$ does not contain any
copies of $G_2$, we see that $r(G_1,G_2)\ge 1+2(d_2-1)=2d_2-1 \ge
d_1+d_2$. This proves (ii).
\par Finally we consider (iii).
Suppose that $G_1$ is a connected graph of order $m$ and $d_2>m$. By
Lemma 2.2, we may construct a graph
$$G=\cases K_{m-1}\cup H_1&\t{if $2\mid
(d_1+d_2-m)$,}\\K_{m-2}\cup H_2&\t{if $2\nmid
(d_1+d_2-m)$,}\endcases$$ where $H_1$ is a $d_1-1$-regular graph of
order $d_1+d_2-m$ and $H_2$ is a $d_1-1$-regular graph of order
$d_1+d_2-m+1$. It is easily seen that $G$ does not contain any
copies of $G_1$ and
$$\Delta(\overline G)=\cases d_2-1&\t{if $2\mid (d_1+d_2-m)$ or
$d_1\not=m-1$,}
\\d_2&\t{if $2\nmid (d_1+d_2-m)$ and $d_1=m-1$.}\endcases$$
If $2\mid (d_1+d_2-m)$ or $d_1\not=m-1$, then $\overline G$ does not
contain any copies of $G_2$ and so $r(G_1,G_2)\ge 1+|V(G)|=d_1+d_2$.
Now assume $2\nmid (d_1+d_2-m)$ and $d_1=m-1$. For $v_0\in V(H_2)$
we have $d_{\overline G}(v_0)=d_2-1$. Suppose that
$v_1,\ldots,v_{m-2}\in V(G)$ and $v_1,\ldots,v_{m-2}$ induce a copy
of $K_{m-2}$. Then $\{v_1,\ldots,v_{m-2}\}$ is an independent set in
$\overline G$ and $d_{\overline G}(v_i)=d_2$ for $i=1,2,\ldots,m-2$.
If $G_2$ has two vertices $u$ and $v$ such that $d(v)=\Delta(G_2)$
and $d(u,v)=3$, we  see that $\overline G$ does not contain any
copies of $G_2$ and so $r(G_1,G_2)\ge 1+|V(G)|=d_1+d_2$. This proves
(iii) and the lemma is proved.

 \pro{Lemma 2.4 ([10, Theorem 2.1])} Let $p,n\in\Bbb N$ with $p\geq n-1\geq 1$. Then
$ex(p;K_{1,n-1})=\lfloor\f{(n-2)p}2\rfloor$.
\endpro

\pro{Lemma 2.5 ([10, Theorem 3.1])} Let $p,n\in\Bbb N$ with $p\geq
n\geq 5$. Let $r\in\{ 0,1,\ldots,n-2\}$ be given by $p\e
r\mod{n-1}$. Then
$$ex(p;T_n')= \cases \big\lfloor\f{(n-2)(p-1)-r-1}2\big\rfloor &\t{if}\
n\ge 7\ \t{and}\
2\le r\le n-4,\\
\f{(n-2)p-r(n-1-r)}2&\t{otherwise.}\endcases$$
\endpro
\pro{Lemma 2.6 ([10, Theorems 4.1-4.3])} Let $p,n\in\Bbb N$ with
$p\ge n\ge 6$, and let $p=k(n-1)+r$ with $k\in\Bbb N$ and
$r\in\{0,1,n-5,n-4,n-3,n-2\}.$ Then $$\aligned ex(p;T_n^*)&=\cases
\f{(n-2)(p-2)}2+1&\t{if $n>6$ and $r=n-5$,}
\\\f{(n-2)p-r(n-1-r)}2&\t{otherwise.}
\endcases\endaligned$$\endpro
\pro{Lemma 2.7 ([10, Theorem 4.4])} Let $p,n\in\Bbb N$, $p\ge n\ge
11$, $r\in\{2,3,\ldots,n-6\}$ and $p\e r\mod{n-1}$. Let
$t\in\{0,1,\ldots,r+1\}$ be given by $n-3\e t\mod{r+2}$. Then
$$\aligned ex(p;T_n^*)=\cases \lfloor\f{(n-2)(p-1)-2r-t-3}2
\rfloor&\t{if $r\ge 4$ and $2\le t\le
r-1$,}\\\f{(n-2)(p-1)-t(r+2-t)-r-1}2
&\t{otherwise}.\endcases\endaligned$$\endpro
 \pro{Lemma 2.8 ([10, Theorem 4.5])}
Let $p,n\in\Bbb N$ with $6\le n\le 10$ and $p\ge n$, and let
$r\in\{0,1,\ldots,n-2\}$ be given by $p\e r\mod{n-1}$.
\par $(\t{\rm i})$ If $n=6,7$, then
$ex(p;T_n^*)=\f{(n-2)p-r(n-1-r)}2$.
\par $(\t{\rm ii})$ If $n=8,9$, then
$$ex(p;T_n^*)=\cases \f{(n-2)p-r(n-1-r)}2&\t{if $r\not=n-5$,}
\\\f{(n-2)(p-2)}2+1&\t{if $r=n-5$.}\endcases$$
\par $(\t{\rm iii})$ If $n=10$, then
$$ex(p;T_n^*)=\cases 4p-\f{r(9-r)}2&\t{if $r\not=4,5$,}
\\4p-7&\t{if $r=5$,}\\4p-9&\t{if $r=4$.}
\endcases$$
\endpro
 \pro{Lemma 2.9} Let $p,m\in\Bbb N$ with $p\ge m\ge 5$,
and $T_m\in\{P_m,K_{1,m-1},$ $T_m',T_m^*\}$. Then $ex(p;T_m)\le
\f{(m-2)p}2$. Moreover, if $m-1\nmid p$ and
$T_m\in\{P_m,T_m',T_m^*\}$, then $ex(p;T_m)\le \f{(m-2)(p-1)}2$.
\endpro
Proof. This is immediate from (1.1) and Lemmas
2.4-2.8.

\pro{Lemma 2.10} Let $m,n\in\Bbb N$ with $m,n\ge 5$. Let $G_m$ be a
connected graph on $m$ vertices. If $m+n-5=(m-1)x+(m-2)y$ for some
nonnegative integers $x$ and $y$, then $r(G_m,T_n)\ge m+n-4$ for
$T_n\in\{K_{1,n-1},T_n',T_n^*\}$.
\endpro Proof.
Let $G=xK_{m-1}\cup yK_{m-2}$. Then $|V(G)|=m+n-5$, $\Delta(G)\le
m-1$ and $\Delta(\overline G)\le n-3$. Clearly, $G$ does not contain
$G_m$ as a subgraph, and $\overline G$ does not contain $T_n$ as a
subgraph. So the result is true.
 \pro{Lemma 2.11 ([7, Theorem 8.3, pp.11-12])} Let $a,b,n\in\Bbb N$. If $a$ is coprime to $b$ and $n\ge (a-1)(b-1)$,
then there are two nonnegative integers $x$ and $y$ such that
$n=ax+by$.\endpro
 \pro{Conjecture 2.12} Let $p,n\in\Bbb N$, $p\ge
n\ge 5$, $p=k(n-1)+r$, $k\in\Bbb N$ and $r\in\{0,1,\ldots,n-2\}$.
Let $T_n\not=K_{1,n-1},T_n'$ be a tree on $n$ vertices. Then
$ex(p;T_n)\le ex(p;T_n^*)$. Hence:
\par $(\t{\rm
i})$ if $r\in\{0,1,n-4,n-3,n-2\}$, then
$$ex(p;T_n)=\f{(n-2)p-r(n-1-r)}2.$$
\par $(\t{\rm
ii})$ if $2\le r\le n-5$, then
$$ ex(p;T_n)\le\f{(n-2)(p-1)-r-1}2.$$
\endpro
\par We note that
$$ex(p;T_n)\ge e(kK_{n-1}\cup K_r)=\f{(n-2)p-r(n-1-r)}2=ex(p;P_n).
$$

\pro{Definition 2.13} For $n\ge 5$ let $T_n$ be a tree on $n$
vertices. View $T_n$ as a bipartite graph with $s_1$ and $s_2$
vertices in the bipartition. Define $\alpha_2(T_n)=\t{max}\
\{s_1,s_2\}$.\endpro

 \pro{Conjecture 2.14} Let $p,n\in\Bbb N$ with $p\ge
n\ge 5$. Let $T_n^{(1)}$ and $T_n^{(2)}$ be two trees on $n$
vertices. If $\alpha_2(T_n^{(1)})< \alpha_2(T_n^{(2)})$, then
$ex(p;T_n^{(1)})\le ex(p;T_n^{(2)})$.
\endpro

 \section*{3. The Ramsey number $r(G_n,T_n^*)$}\par\q

 \pro{Lemma 3.1} Let $n\in\Bbb N$, $n\ge
6$, and let $G_n$ be a connected graph on $n$ vertices such that
$ex(2n-5;G_n)<n^2-5n+4.$ Then $r(G_n,T_n^*)=2n-5.$\endpro Proof. As
$2K_{n-3}$ does not contain any copies of $G_n$ and
$\overline{2K_{n-3}}=K_{n-3,n-3}$ does not contain any copies of
$T_n^*$, we see that $r(G_n,T_n^*)>2(n-3).$ By Lemma 2.6 we have
$$ex(2n-5;T_n^*)=\f{(n-2)(2n-5)-3(n-4)}2=n^2-6n+11.$$ Thus,
$$\aligned ex(2n-5;G_n)+ex(2n-5;T_n^*)&<n^2-5n+4+n^2-6n+11
\\&=2n^2-11n+15=\b{2n-5}2.\endaligned$$
Appealing to Lemma 2.1 we obtain $r(G_n,T_n^*)\le 2n-5$. So
$r(G_n,T_n^*)=2n-5$ as asserted.
 \pro{Theorem 3.2} Let $n\in\Bbb N$
with $n\ge 8$. Then
$$r(P_n,T_n^*)=r(T_n',T_n^*)=r(T_n^*,T_n^*)=2n-5.$$\endpro
Proof. By Lemma 2.6,
$$ex(2n-5;T_n^*)=\f{(n-2)(2n-5)-3(n-4)}2=n^2-6n+11<n^2-5n+4.$$ By
Lemma 2.5,
$$\aligned ex(2n-5;T_n')&=\Big\lfloor\f{(n-2)(2n-6)-(n-4)-1}2\Big\rfloor
=\Big\lfloor n^2-\f{11}2n+\f{15}2\Big\rfloor\\&\le
n^2-\f{11}2n+\f{15}2 <n^2-5n+4.\endaligned$$ By (1.1),
$$ex(2n-5;P_n)=\b{n-1}2+\b{n-4}2=n^2-6n+11<n^2-5n+4.$$ Thus applying
Lemma 3.1 we deduce the result.

 \pro{Conjecture 3.3} Let $n\in\Bbb N$, $n\ge
8$, and let $T_n\neq K_{1,n-1}$ be a tree on $n$ vertices. Then
$r(T_n,T_n^*)=2n-5.$\endpro
 \par\q\par {\bf Remark 3.4} Let
$n\in\Bbb N$ with $n\ge 4.$ From [6, Theorem 3.1(ii)] we know that
$r(K_{1,n-1},T_n^*)=2n-3.$

\section*{4. The Ramsey number $r(G_m,T_n^*)$ for $m<n$}\par\q
\pro{Theorem 4.1} Let $m,n\in\Bbb N$, $n>m\ge 5$ and $m-1\mid n-3$.
Let $G_m$ be a connected graph of order $m$ such that
$ex(m+n-3;G_m)\le\f{(m-2)(m+n-3)}2$ or
$G_m\in\{P_m,K_{1,m-1},T_m',T_m^*\}$. Then
$r(G_m,T_n^*)=m+n-3.$\endpro
 Proof. By Lemma 2.9 we may assume that
$ex(m+n-3;G_m)\le\f{(m-2)(m+n-3)}2.$ Suppose that $n-3=k(m-1).$
Clearly $(k+1)K_{m-1}$ does not contain $G_m$ as a subgraph and
$\overline{(k+1)K_{m-1}}$ does not contain $T_n^*$ as a subgraph.
Thus
$$r(G_m,T_n^*)>(k+1)(m-1)=m+n-4.$$ Since $1\le m-4\le n-6$,
using Lemma 2.9 we see that
$$\aligned
ex(m+n-3;T_n^*)\le \f{(n-2)(m+n-4)}2.\endaligned$$ Thus,
$$\aligned
&ex(m+n-3;G_m)+ex(m+n-3;T_n^*)\\&\le\f{(m-2)(m+n-3)}2+\f{(n-2)(m+n-4)}2
\\&<\f{(m-2+n-2)(m+n-3)}2=\b{m+n-3}2.\endaligned$$ Hence, by Lemma
2.1, $r(G_m,T_n^*)\le m+n-3$, and the result follows.
 \pro{Lemma
4.2} Let $m,n\in\Bbb N,n>m\ge 7$ and $m-1\nmid n-3.$ Let $G_m$ be a
connected graph of order $m$ such that
$ex(m+n-4;G_m)\le\f{(m-2)(m+n-4)}2$ or
$G_m\in\{P_m,K_{1,m-1},T_m',T_m^*\}$. Then $r(G_m,T_n^*)\le
m+n-4.$\endpro
 Proof. By Lemma 2.9, we may assume that $ex(m+n-4;G_m)\le
\f{(m-2)(m+n-4)}2.$ As $m+n-4=n-1+m-3$ and $m-1\nmid (n-3),$ we see
that $2\le m-3\le n-4$ and $m-3\neq n-5.$ Thus, applying Lemmas 2.6-
2.8,
$$ex(m+n-4;T_n^*)<\f{(n-3)(m+n-4)}2.$$
Hence, $$\aligned
&ex(m+n-4;G_m)+ex(m+n-4;T_n^*)\\&<\f{(m-2)(m+n-4)}2+\f{(n-3)(m+n-4)}2
=\b{m+n-4}2.\endaligned $$ Applying Lemma 2.1, we obtain the result.
 \pro{Theorem 4.3} Let $m,n\in\Bbb N,n>m\ge 7$ and $m-1\nmid
(n-3).$ Let $G_m$ be a connected graph of order $m$ such that
$ex(m+n-4;G_m)\le\f{(m-2)(m+n-4)}2$ or $G_m\in\{P_m,T_m',T_m^*\}$.
If $m+n-5=(m-1)x+(m-2)y$ for some $x,y\in\{0,1,2,\ldots\}$, then
$r(G_m,T_n^*)=m+n-4.$\endpro Proof. By Lemma 4.2, $r(G_m,T_n^*)\le
m+n-4$, and by Lemma 2.10, $r(G_m,T_n^*)\ge m+n-4$. Thus the result
follows.

\pro{Theorem 4.4} Suppose $m,n\in\Bbb N$, $n>m\ge 7$,
$n=k(m-1)+b=q(m-2)+a$, $k,q\in\Bbb N$, $a\in\{0,1,\ldots,m-3\}$ and
$b\in\{0,1,\ldots,m-2\}-\{3\}$. Let $G_m$ be a connected graph of
order $m$ such that $ex(m+n-4;G_m)\le\f{(m-2)(m+n-4)}2$ or
$G_m\in\{P_m,T_m',T_m^*\}$. If one of the conditions:
$$\align &(\t{\rm i})\q b\in\{1,2,4\},
\\ &(\t{\rm ii})\q b=0\qtq{and}k\ge 3,
\\&(\t{\rm iii})\q n\ge (m-3)^2+2,
\\&(\t{\rm iv})\q n\ge m^2-1-b(m-2),
\\ &(\t{\rm v})\q\ a\ge 3\qtq{and}n\ge (a-4)(m-1)+4\qq\qq\qq\qq\qq
\endalign$$
holds, then $r(G_m,T_n^*)=m+n-4$.\endpro Proof. For $b\in\{1,2,4\}$,
$$m+n-5=\cases (k-2)(m-1)+3(m-2) &\t{if}\ b=1,
\\(k-1)(m-1)+2(m-2) &\t{if}\ b=2,
\\(k+1)(m-1) &\t{if}\ b=4.\endcases$$
For $b=0$ and $k\ge 3$ we have $m+n-5=(k-3)(m-1)+4(m-2)$. For $n\ge
(m-3)^2+2$,  we have $m+n-5\ge (m-2)(m-3)$ and so
$m+n-5=(m-1)x+(m-2)y$ for some $x,y\in\{0,1,2,\ldots\}$ by Lemma
2.11. For $n\ge m^2-1-b(m-2)$ we have $k\ge m+1-b$ and
$m+n-5=(k+b-m-1)(m-1)+(m+3-b)(m-2)$. For $a\ge 3$ and $n\ge
(a-4)(m-1)+4$ we have $q\ge a-4$ and
$m+n-5=(a-3)(m-1)+(q+4-a)(m-2)$. Combining all the above with
Theorem 4.3, we obtain the result.

 \pro{Theorem
4.5} Suppose that $m,n\in\Bbb N$, $n>m\ge 7$ and $m-1\nmid n-3$.
  Then
$$\aligned &r(K_{1,m-1},T_n^*)=m+n-4,
\\&r(T_m',T_n^*)=
m+n-4\ \t{or}\ m+n-5,
\\& m+n-6\le r(T_m^*,T_n^*)\le m+n-4.\endaligned$$
\endpro
 Proof. From Lemma 4.2, $r(T_m,T_n^*)\le m+n-4$ for
 $T_m\in\{K_{1,m-1},T_m',T_m^*\}$. By Lemma 2.3,
$r(K_{1,m-1},T_n^*)\ge m-1+n-3$, $r(T_m',T_n^*)\ge m-2+n-3$
$(n>m+1)$ and $r(T_m^*,T_n^*)\ge m-3+n-3.$ By Theorem 4.4,
$r(T_m',T_n^*)=m+n-4$ for $n=m+1,m+3$.
 Thus the theorem is
proved.

\pro{Theorem 4.6} Suppose that $m,n\in\Bbb N$, $n>m\ge 7$,
$n=k(m-1)+b$, $k\in\Bbb N$, $b\in\{0,1,\ldots,m-2\}$, $b\neq 3$ and
$\f{m-b}2\le k\le m+2-b$. Let $G_m$ be a connected graph
 of order  $m$ such that $ex(m+n-4;G_m)\le \f 12(m-2)(m+n-4)$ or
$G_m\in\{P_m,T_m^*\}$. Then $r(G_m,T_n^*)=m+n-4$ or $m+n-5$.\endpro
Proof. By Lemma 4.2 we only need to show that $r(G_m,T_n^*)>m+n-6$.
Set $G=(2k+b-m)K_{m-2}\cup (m+2-b-k)K_{m-3}$. Then
$|V(G)|=(2k+b-m)(m-2)+(m+2-b-k)(m-3)=m+n-6$. We also have
$\Delta(G)\le m-2$ and $\Delta(\overline G)\le m+n-6-(m-3)=n-3$. Now
it is clear that $G_m$ is not a subgraph of $G$ and that $T_n^*$ is
not a subgraph of $\overline G$. So $r(G_m,T_n^*)>|V(G)|$, which
completes the proof.
\par\q
\newline\par{\bf Remark 4.7} If $p\ge m\ge 6$ and $T_m$ is a
tree on $m$ vertices with a vertex adjacent to at least
$\lfloor\f{m-1}2\rfloor$ vertices of degree $1$, in [9] Sidorenko
proved that $ex(p;T_m)\le \f{(m-2)p}2$. Thus, $G_m$ can be replaced
with $T_m$ in Lemma 4.2, Theorems 4.1, 4.3, 4.4 and 4.6.

\section*{5. The Ramsey number $r(G_m,T_n')$ for $m<n$}\par\q
\pro{Theorem 5.1} Let $m,n\in\Bbb N$, $n>m\ge 6$ and $m-1\mid n-3.$
Suppose that $G_m$ is a connected graph of order $m$ satisfying
$ex(m+n-3;G_m)\le\f{(m-2)(m+n-3)+m+n-4}2$ or $G_m\in \{T_m^*,P_m\}$.
Then $r(G_m,T_n')=m+n-3.$\endpro
 Proof. By Lemma 2.9 we may assume that
 $$ex(m+n-3;G_m)\le(m-2)(m+n-3)/2+(m+n-4)/2.$$
 Suppose $n-3=k(m-1)$ and
$G=(k+1)K_{m-1}.$ Then $|V(G)|=m+n-4$ and $\Delta(\overline G)=n-3.$
Clearly, $G_m$ is not a subgraph of $G$ and $T_n'$ is not a subgraph
of $\overline G.$ Thus $r(G_m,T_n')>m+n-4.$
 Since
$m-1\mid (n-3),$ we have $n\ge m+2$ and so $4\le m-2\le n-4$. Hence,
using Lemma 2.5, $ex(m+n-3;T_n')=\lfloor
\f{(n-2)(m+n-4)-(m-1)}2\rfloor<\f{(n-2)(m+n-3)-(m+n-4)}2.$ Therefore
$$ex(m+n-3;G_m)+ex(m+n-3;T_n')<\b{m+n-3}2.$$
Applying Lemma 2.1, we see that $r(G_m,T_n')\le m+n-3$, so the
result follows.
 \pro{Lemma 5.2} Let $m,n\in\Bbb N,n>m\ge 6$ and
$m-1\nmid n-3.$ Suppose that $G_m$ is a connected graph of order $m$
satisfying $ex(m+n-4;G_m)<\f{(m-2)(m+n-4)}2$ or $G_m\in
\{T_m^*,P_m\}$. Then $r(G_m,T_n')\le m+n-4.$\endpro
 Proof.  Since
$m-1\nmid n-3$, $m-1\nmid m+n-4.$ Thus, applying Lemma 2.9,
 $ex(m+n-4;T_m^*)\le (m-2)(m+n-5)/2$ and
  $ex(m+n-4;P_m)\le(m-2)(m+n-5)/2$.
  As $n>m$, $3\le m-3\le n-4$. By Lemma 2.5,
 $ex(m+n-4;T_n')=\lfloor\f{(n-2)(m+n-5)-(m-2)}2\rfloor
 \le \f{(n-2)(m+n-5)-(m-2)}2$. Thus
 $$ex(m+n-4;G_m)+ex(m+n-4;T_n')<\f{(m-2+n-2)(m+n-5)}2=\b{m+n-4}2.$$
 This, together with Lemma 2.1, yields the result.

 \pro{Theorem 5.3} Let $m,n\in\Bbb N$, $n>m\ge 6$ and
$m-1\nmid (n-3).$  Then $r(T_m^*,T_{m+1}')=2m-3$ and
$r(T_m^*,T_n')=m+n-4$ or $m+n-5$ for $n\ge m+3$. Suppose that $G_m$
is a connected graph of order $m$ satisfying
$ex(m+n-4;G_m)<\f{(m-2)(m+n-4)}2$ or $G_m\in \{T_m^*,P_m\}$. If
$m+n-5=(m-1)x+(m-2)y$ for some nonnegative integers $x$ and $y$,
then $r(G_m,T_n')=m+n-4.$\endpro
 Proof. By Lemma 2.3, $r(T_m^*,T_{m+1}')\ge 2(m-1)-1=2m-3$ and
 $r(T_m^*,T_n')\ge m-3+n-2$ for $n\ge m+3$. By Lemma 5.2,
 $r(G_m,T_n')\le m+n-4$. Thus, $r(T_m^*,T_{m+1}')=2m-3$.
 Applying Lemma 2.10 we deduce
 the remaining result.
\par\q
\par From Theorem 5.3 and the proof of Theorem 4.4 we deduce the
following result.
 \pro{Theorem 5.4} Suppose $m,n\in\Bbb N$, $n>m\ge
6$, $n=k(m-1)+b=q(m-2)+a$, $k,q\in\Bbb N$, $a\in\{0,1,\ldots,m-3\}$
and $b\in\{0,1,\ldots,m-2\}-\{3\}$. Let $G_m$ be a connected graph
of order $m$ such that $ex(m+n-4;G_m)<\f{(m-2)(m+n-4)}2$ or
$G_m\in\{P_m,T_m^*\}$. If one of the  conditions:
$$\align &(\t{\rm i})\q b\in\{1,2,4\},
\\ &(\t{\rm ii})\q b=0\qtq{and}k\ge 3,
\\&(\t{\rm iii})\q n\ge (m-3)^2+2,
\\&(\t{\rm iv})\q n\ge m^2-1-b(m-2),
\\ &(\t{\rm v})\q\ a\ge 3\qtq{and}n\ge (a-4)(m-1)+4\qq\qq\qq\qq\qq
\endalign$$
holds, then $r(G_m,T_n')=m+n-4$.\endpro

\section*{6. The Ramsey number $r(T_m,K_{1,n-1})$ for $m<n$}\par\q
The following two propositions are known.
 \pro{Proposition 6.1
([1])} Let $m,n\in\Bbb N$ with $m\ge 3$ and $m-1\mid n-2.$ Let $T_m$
be a tree on $m$ vertices. Then $r(T_m,K_{1,n-1})=m+n-2.$\endpro
 \pro{Proposition 6.2 ([6, Theorem 3.1])} Let $m,n\in\Bbb N, m\ge
3$ and $n=k(m-1)+b$ with $k\in\Bbb N$ and
$b\in\{0,1,\ldots,m-2\}-\{2\}.$ Let $T_m\not=K_{1,m-1}$ be a tree on
$m$ vertices. Then $r(T_m,K_{1,n-1})\le m+n-3$. Moreover, if $k\ge
m-b$, then $r(T_m,K_{1,n-1})=m+n-3.$\endpro

\pro{Theorem 6.3} Let $m,n\in\Bbb N$, $n\ge m\ge 3$, $m-1\nmid
(n-2)$, $n=q(m-2)+a$, $q\in\Bbb N$ and $a\in\{2,3,\ldots,m-3\}.$ Let
$T_m\not=K_{1,m-1}$ be a tree on $m$ vertices. If $n\ge
(a-3)(m-1)+3,$ then $r(T_m,K_{1,n-1})=m+n-3.$\endpro
 Proof. Since
$q(m-2)=n-a\ge(a-3)(m-2)$ we have $q\ge a-3$.
 Set $G=(a-2)K_{m-1}\cup (q-(a-3))K_{m-2}$. Then
 $|V(G)|=(a-2)(m-1)+(q-(a-3))(m-2)=m+n-4$ and $\Delta(\overline G)\le
 n-2.$ Clearly, $T_m$ is not a subgraph of $G$
and $K_{1,n-1}$ is not a subgraph of $\overline G$. Thus
$r(T_m,K_{1,n-1})>|V(G)|=m+n-4.$ By Proposition 6.2,
$r(T_m,K_{1,n-1})\le m+n-3.$ So $r(T_m,K_{1,n-1})=m+n-3.$ This
proves the theorem.

\pro{Theorem 6.4} Let $m,n\in\Bbb N$ with $n>m\ge 5$ and $m-1\nmid
(n-2)$. Then $r(T_m^*,K_{1,n-1}) =m+n-3$ or $m+n-4$. Moreover, if
$m+n-4=(m-1)x+(m-2)y+2(m-3)z$ for some nonnegative integers $x,y$
and $z$, then $r(T_m^*,K_{1,n-1})=m+n-3$.
\endpro Proof. By Proposition 6.2, $r(T_m^*,K_{1,n-1})\le m+n-3$.
By Lemma 2.3 we have $r(T_m^*,K_{1,n-1})\ge m+n-4$.  If
$m+n-4=(m-1)x+(m-2)y+2(m-3)z$ for some nonnegative integers $x,y$
and $z$, setting $G=xK_{m-1}\cup yK_{m-2}\cup zK_{m-3,m-3}$ we find
$\Delta(\overline G)\le n-2$. Clearly, $G$ does not contain any
copies of $T_m^*$, and $\overline G$ does not contain any copies of
$K_{1,n-1}$. Thus, $r(T_m^*,K_{1,n-1})>|V(G)|=m+n-4$ and so
$r(T_m^*,K_{1,n-1})=m+n-3$. This proves the theorem.

\enddocument
\bye